\documentclass[12pt,english]{article}
\usepackage{babel}
\usepackage[cp1250]{inputenc}
\usepackage[T1]{fontenc}
\usepackage{amsfonts}
\usepackage{amsmath}
\usepackage{amsthm}
\usepackage{hyperref}
\newcommand{\no}{\noindent}
\newcommand{\rr}{\mathbb{R}}
\newcommand{\hd}{\hspace{0.2cm}}
\newcommand{\hdd}{\hspace{0.5cm}}
\newcommand{\lap }{\Delta}

\newcommand{\dx}{\hspace{0.1cm} dx}

\newcommand{\dt}{\hspace{0.1cm} dt}
\newcommand{\ur}{u_{\rho}}
\newcommand{\urm}{u^{-}_{\rho}}
\newcommand{\uh}{u_{\theta}}
\newcommand{\uz}{u_{z}}

\newcommand{\urr}{\poch{\ur}{\rho}}
\newcommand{\urz}{\poch{\ur}{z}}

\newcommand{\uhr}{\poch{\uh}{\rho}}
\newcommand{\uhz}{\poch{\uh}{z}}

\newcommand{\uzr}{\poch{\uz}{\rho}}
\newcommand{\uzz}{\poch{\uz}{z}}

\newcommand{\hr}{h_{\rho}}
\newcommand{\hh}{h_{\theta}}
\newcommand{\hz}{h_{z}}
\newcommand{\oo}{\omega}
\newcommand{\oor}{\oo_{\rho}}
\newcommand{\oh}{\oo_{\theta}}
\newcommand{\oz}{\oo_{z}}
\newcommand{\ort}{\poch{\oor}{t}}
\newcommand{\orr}{\poch{\oor}{\rho}}
\newcommand{\orz}{\poch{\oor}{z}}

\newcommand{\oht}{\poch{\oh}{t}}
\newcommand{\ohr}{\poch{\oh}{\rho}}
\newcommand{\ohz}{\poch{\oh}{z}}

\newcommand{\ozt}{\poch{\oz}{t}}
\newcommand{\ozr}{\poch{\oz}{\rho}}
\newcommand{\ozz}{\poch{\oz}{z}}

\newcommand{\grr}{g_{\rho}}
\newcommand{\gh}{g_{\theta}}
\newcommand{\gz}{g_{z}}
\newcommand{\uhq}{\uh^{q}}
\newcommand{\uhqj}{\uh^{q-1}}
\newcommand{\uhqd}{\uh^{q-2}}
\newcommand{\uhrk}{(\uhr)^{2}}
\newcommand{\uhzk}{(\uhz)^{2}}
\newcommand{\ib}{\intl{}{B_{2}}}
\newcommand{\qo}{\frac{1}{q}}
\newcommand{\pt}{\frac{d}{dt}}
\newcommand{\uhqr}{\poch{\uhq}{\rho}}
\newcommand{\uhqz}{\poch{\uhq}{z}}
\newcommand{\ibt}{\intl{}{\widetilde{B}_{2}}}

\newcommand{\uhqprk}{(\poch{\uh^{q/2 }}{\rho})^{2}}
\newcommand{\uhqpzk}{(\poch{\uh^{q/2 }}{z})^{2}}

\newcommand{\uhqpjk}{(\uh^{q/2 -1})^{2}}
\newcommand{\pqp}{\big[\uhqprk + \uhqpzk\big]}
\newcommand{\hhq}{\hh^{q}}
\newcommand{\epj}{\varepsilon_{1}}
\newcommand{\ep}{\varepsilon}
\newcommand{\epd}{\varepsilon_{2}}
\newcommand{\ids}{\intl{t}{t_{0}-\tau}d(s)ds}
\newcommand{\rde}{\rho^{2-\ep}}
\newcommand{\rte}{\rho^{3-\ep}}
\newcommand{\rje}{\rho^{1-\ep}}
\newcommand{\rce}{\rho^{4-\ep}}
\newcommand{\ohkrde}{\frac{\oh^{2}}{\rde}}
\newcommand{\ohkrd}{\frac{\oh^{2}}{\rho^{2}}}
\newcommand{\pohkrde}{\frac{1}{2} \frac{d}{dt} \ib \ohkrde}
\newcommand{\pohkrd}{\frac{1}{2} \frac{d}{dt} \ib \ohkrd}
\newcommand{\jd}{\frac{1}{2}}
\newcommand{\ohkr}{\poch{\oh^{2}}{\rho}}
\newcommand{\ohkz}{\poch{\oh^{2}}{z}}
\newcommand{\ohk}{\oh^{2}}
\newcommand{\ohrk}{(\ohr)^{2}}
\newcommand{\ohzk}{(\ohz)^{2}}
\newcommand{\epp}{\frac{\ep}{2}}

\newcommand{\uhpr}{\frac{\uh}{\rho}}
\newcommand{\uhkz}{\poch{\uh^{2}}{z}}
\newcommand{\uhk}{\uh^{2}}
\newcommand{\ohrjerk}{\big[ \poch{}{\rho} \big(\frac{\oh}{\rje}\big) \big]^{2}}
\newcommand{\ohrjezk}{\big[ \poch{}{z} \big(\frac{\oh}{\rje}\big) \big]^{2}}
\newcommand{\ohrjrk}{\big[ \poch{}{\rho} \big(\frac{\oh}{\rho}\big) \big]^{2}}
\newcommand{\ohrjzk}{\big[ \poch{}{z} \big(\frac{\oh}{\rho}\big) \big]^{2}}
\newcommand{\roe}{\frac{1}{\rho^{\ep}}}
\newcommand{\rdeo}{\frac{1}{\rde}}
\newcommand{\uhtrd}{\frac{\uh^{3}}{\rho^{2}}}
\newcommand{\jc}{\frac{1}{4}}
\newcommand{\puhcrd}{\frac{d}{dt} \ib \frac{\uh^{4}}{\rho^{2}}}
\newcommand{\uhc}{\uh^{4}}
\newcommand{\uhcrd}{\frac{\uhc}{\rho^{2}}}
\newcommand{\uhcrc}{\frac{\uhc}{\rho^{4}}}
\newcommand{\uhdrd}{\frac{\uhk}{\rho^{2}}}
\newcommand{\uhcurrt}{\frac{\uhc \ur}{\rho^{3}}}
\newcommand{\uhcurbrt}{\frac{\uhc |\ur|}{\rho^{3}}}
\newcommand{\rod}{\frac{1}{\rho^{2}}}
\newcommand{\uhkrzk}{ \Big[\poch{}{z} \big( \frac{\uhk}{\rho} \big) \Big]^{2}}
\newcommand{\uhkrrk}{ \Big[\poch{}{\rho} \big( \frac{\uhk}{\rho} \big) \Big]^{2}}

\newcommand{\eqq}[2]{\begin{equation} #1 \label{#2} \end{equation}}
\newcommand{\m}[1]{\mbox{#1 }}

\newcommand{\poch}[2]{\frac{\partial #1 }{\partial #2 }}
\newcommand{\pochd}[2]{\frac{\partial^{2} #1 }{\partial {#2}^{2} }}
\newcommand{\intl}[2]{\int\limits^{#1}_{#2}}

\newcommand{\bl}[1]{\mathbf{#1}}
\newcommand{\lqq}[1]{\| #1 \|_{q} }
\newcommand{\lqqq}[1]{\| #1 \|^{q}_{q} }

\newcommand*{\divv}{\mathop{\mathrm{div}}\limits}
\newcommand*{\curl}{\mathop{\mathrm{curl}}\limits}

\newtheorem{remark}{\textbf{Remark}}[section]

\title{Remarks on Interior Regularity Criterion for an Axially Symmetric
Suitable Weak Solution to the Navier Stokes Equations}

\author{Adam Kubica}

\date{January 3, 2010}

\newtheorem{common}{{\textbf {Theorem}}}

\begin{document}

\maketitle

\begin{center}
Faculty  of Mathematics and Information Science \\
Warsaw University of Technology \\
Pl. Politechniki 1, Warsaw 00-661 \\
A.Kubica@mini.pw.edu.pl \\
\end{center}

\abstract{ We show that if v is an axially symmetric suitable weak
solution to the Navier Stokes equations (in the sense of L.
Caffarelli, R. Kohn $\&$ L. Nirenberg) such that the radial
component of v  has a higher regularity (i.e. satisfies weighted
Serrin-Prodi type condition), then all components of v are regular.}

\vspace{2cm}

\no \textbf{Introduction. } In paper \cite{pokornynestupa} is proved
a result concerning conditional regularity of an axially symmetric
suitable weak solutions of Navier-Stokes equations. The authors show
that if $v_{r}$ is the  radial component of the velocity and
satisfies Serrin-Prodi type condition, i.e.

\[
\int^{T}_{0} \Big( \int_{\Omega}
|v_{r}|^{a}dx\Big)^{\frac{b}{a}}dt<\infty ,
\]
\no with $\frac{3}{a}+\frac{2}{b} \leq 1$, $a\in(3,\infty]$, $b\in
[2, \infty]$, then $v$ is \textit{regular}. In this paper we modify
  the proof and obtain the same result under a more general
  assumption: weighted Serrin-Prodi type condition.

We suppose that $\Omega$ is either $\rr^{3}$ or an axially symmetric
bounded domain with smooth boundary and denote $Q_{T}=\Omega\times
(0,T)$ for $T>0$.  We consider the following problem

\eqq{ \poch{\bl{v}}{t}+(\bl{v} \cdot \nabla v )\bl{v}= \bl{f}-
\nabla  p + \nu \lap \bl{v} \m{ \/ in \/ } Q_{T} }{ta}
\eqq{\divv{\bl{u}}=0 \m{ \/ in \/ } Q_{T} }{tb} \eqq{\bl{v}=0 \m{ \/
on \/ } \partial \Omega \times (0,T)}{tc}
\eqq{\bl{v}_{t=0}=\bl{v_{0}}.}{td}

\no We will further suppose for simplicity that $\bl{f}=0$.
Proceeding similarly as in \cite{pokornynestupa} we can reduce the
above problem to the problem on $B_{2}\times (t_{0}- \tau, t_{0})$.
Then we have $\bl{u}$ which satisfies in a classical sense the
equations

\eqq{ \poch{\bl{u}}{t}+(\bl{u} \cdot \nabla u )\bl{u}= h- \nabla
(\eta p )+ \nu \lap \bl{u} }{a} \eqq{\divv{\bl{u}}=0}{b}

\no  in $B_{2}\times (t_{0}- \tau, t_{0})$. Our goal is to proved
that $\bl{u}$ do not blow up at $t=t_{0}$. Therefore we have to
prove appropriate estimates for solution $\bl{u}$ under the
assumption that $\bl{u} $ is axially symmetric and $\ur$ (the radial
component of the velocity) has a higher regularity. It is convenient
to write the equation (\ref{a})  in cylindrical coordinates

\eqq{ \poch{\ur}{t} +\ur \urr + \uz \urz - \frac{1}{\rho} \uh^{2}+
\poch{(\eta p) }{\rho}= \hr+ \nu \big[\frac{1}{\rho} \poch{}{\rho}
(\rho \urr) + \pochd{\ur}{z}- \frac{\ur}{\rho^{2}}\big] }{c}

\eqq{ \poch{\uh}{t} +\ur \uhr + \uz \uhz + \frac{1}{\rho} \uh\ur =
\hh+ \nu \big[\frac{1}{\rho} \poch{}{\rho} (\rho \uhr) +
\pochd{\uh}{z}- \frac{\uh}{\rho^{2}}\big] }{d}

\eqq{ \poch{\uz}{t} +\ur \uzr + \uz \uzz + \poch{(\eta p)}{z} = \hz+
\nu \big[\frac{1}{\rho} \poch{}{\rho} (\rho \uzr) + \pochd{\uz}{z}
\big] .}{e}

\no The equation of continuity has the following form in cylindrical
coordinates

\eqq{\urr + \frac{\ur}{\rho}+ \uzz=0.}{f}

\no We put

\eqq{\mbox{\boldmath$\oo$} = \curl{\bl{u}}, \hdd \bl{g}=
\curl{\bl{h}}.}{g}

\no Then we have

\eqq{\oor= - \uhz, \hd \oh= \urz- \uzr, \hd \oz= \frac{1}{\rho}
\poch{(\rho \uh)}{\rho} }{h}

\no Applying operator $\curl$ to equation (\ref{a})we obtain the
system

\eqq{ \ort +\ur \orr + \uz \orz - \urr \oor-\urz \oz = \grr+ \nu
\big[\frac{1}{\rho} \poch{}{\rho} (\rho \orr) + \pochd{\oor}{z}-
\frac{\oor}{\rho^{2}}\big] }{i}

\eqq{ \oht +\ur \ohr + \uz \ohz - \frac{\ur}{\rho} \oh-2
\frac{\uh}{\rho} \oor = \gh+ \nu \big[\frac{1}{\rho} \poch{}{\rho}
(\rho \ohr) + \pochd{\oh}{z}- \frac{\oh}{\rho^{2}}\big] }{j}

\eqq{ \ozt +\ur \ozr + \uz \ozz -\uzr \oor- \uzz \oz= \gz+ \nu
\big[\frac{1}{\rho} \poch{}{\rho} (\rho \ozr) + \pochd{\oz}{z} \big]
.}{k}

Our result is following

\begin{common}
Let $\bl{v}$ be an axially symmetric suitable weak solution to the
problem (\ref{ta})-(\ref{td}) with $\bl{f}=0$. Suppose that there
exists a sub-domain $D$ of $Q_{T}$ such that the radial component
$v_{\varrho}$ of $\bl{v}$ has its negative part $v^{-}_{\varrho}$ in
$L^{b,a}_{\gamma}(D)$ for some $a\in (\frac{3}{2}, \infty]$,  $b\in
(1, \infty)$ such that $\frac{3}{a} +\frac{2}{b}+\gamma \leq 1  $
and $\frac{3}{a} +\frac{2}{b}<2$. Then $\bl{v}$ is regular in $D$.

\end{common}

\vspace{0.3cm}

\no The condition $v^{-}_{\varrho}\in L^{b,a}_{\gamma}(U\times
(0,T))$ means that

\eqq{\intl{T}{0} \Big(\intl{}{U} |\urm \cdot \rho^{\gamma}|^{a} \dx
\Big)^{\frac{b}{a}} \dt <\infty.  }{m}

\vspace{1cm}

\no The prove will be given in several steps.

\no \textbf{Step 1.} Assume that $q$ is even, $t\in (t_{0}-\tau,
t_{0})$ and multiply equation (\ref{d}) by $\uhqj$ and integrate
over $B_{2}$. Then we get

\[
\ib \poch{\uh}{t} \uhqj  +\ib\ur \uhr \uhqj  + \ib\uz \uhz \uhqj +
 \ib\frac{1}{\rho}\ur\uhq
\]
\[
= \ib\hh \uhqj + \nu
\ib\big[\frac{1}{\rho} \uhqj \poch{}{\rho} (\rho \uhr)  + \uhqj
\pochd{\uh}{z}- \frac{\uhq}{\rho^{2}}\big].
\]

\no We have
\[
\ib \poch{\uh}{t} \uhqj= \qo \pt \ib \uhq,
\]
\[
\ib\ur \uhr \uhqj= \qo \ib \ur \uhqr,
\]
\[
\ib\uz \uhz \uhqj= \qo \ib \uz \uhqz,
\]
\[
\ib\frac{1}{\rho} \uhqj \poch{}{\rho} (\rho
\uhr)\overset{\footnotemark}{=} \ibt \uhqj \poch{}{\rho} (\rho \uhr)
\overset{\footnotemark}{=} -\ibt \rho \uhr \poch{\uhqj}{\rho}= (1-q)
\ib \uhrk \uhqd,
\]
\footnotetext{$\widetilde{B}_{2}$ is the ball $B_{2}$ given in
cylindrical coordinates.}
\[
\ib \uhqj \pochd{\uh}{z} \overset{\footnotemark}{=} (1-q) \ib \uhzk
\uhqd.
\]

\no Thus we get

\[
\qo \pt \ib \uhq  +\qo \ib \ur \uhqr  + \qo \ib \uz \uhqz +
 \ib\frac{1}{\rho}\ur\uhq + \nu (q-1) \ib [\uhrk+\uhzk]
\uhqd +\nu\ib \frac{\uhq}{\rho^{2}}
\]
\eqq{
= \ib\hh \uhqj .
}{n}

\no Using (\ref{f}) we get

\[
\ib \ur \uhqr  +  \ib \uz \uhqz = \ibt (\ur  \rho  ) \uhqr  +  \ib  \uz \uhqz= - \ibt (\rho \urr+\ur) \uhq -
  \ib \uzz \uhq
\]
\eqq{
=- \ib ( \urr+\frac{\ur}{\rho} +\uzz) \uhq= 0.
}{o}

\no Clearly we have
\[
[\uhrk+\uhzk] \uhqd= \uhqpjk\uhrk + \uhqpjk \uhzk= (2/q)^{2} \big[\uhqprk + \uhqpzk\big],
\]

\no thus equality (\ref{n}) has the following form

\eqq{
\qo \pt \ib \uhq  +
 \ib\frac{1}{\rho}\ur\uhq + \nu \frac{(q-1)}{(q/2)^{2}} \ib \pqp +\nu\ib \frac{\uhq}{\rho^{2}}
= \ib\hh \uhqj .
}{p}

\no Applying Young inequality\footnote{$ab \leq \qo (\frac{q-1}{q})^{q-1}a^{q}+ b^{q/q-1}$.} we get $\ib\hh \uhqj \leq \qo (\frac{q-1}{q})^{q-1} \ib \hhq + \ib \uhq$. Hence from (\ref{p}) we get the estimate\footnote{$\| \cdot \|_{q}$ denotes $\| \cdot \|_{L^{q}(B_{2})}$.}

\eqq{ \pt  \lqqq{\uh}
 + \nu \frac{4(q-1)}{q}  \ib \pqp +\nu q \ib \frac{\uhq}{\rho^{2}} \leq  q\ib\frac{1}{\rho}\urm\uhq + q\lqqq{\uh} + \lqqq{\hh}.
}{r}

\no Now we shall estimate the first term on the right hand side. We set
\eqq{p= 1+\frac{2a+3b}{2ab-2a-3b}, \hdd s=2\frac{a}{b}+3. }{s}

\no Then $s>3$ and $p>1$, because from the assumption (\ref{l}) we get $3b+ 2a<2ab $. Therefore we may write

\[
\ib\frac{1}{\rho}\urm\uhq = \ib \urm \uh^{q(p-1)/p }\rho^{(2-p)/p} \cdot \uh^{q/p} \rho^{-2/p} \overset{\footnotemark}{\leq } \Big( \ib |\urm|^{\frac{p}{p-1}} \uhq \rho^{\frac{2-p}{p-1}}
 \Big)^{(p-1)/p} \Big(\ib \frac{\uhq}{\rho^{2}} \Big)^{1/p}
\]
\footnotetext{H\"{o}lder: $(\frac{p}{p-1}, p)$.}

\[
=\Big( \epj^{1/(1-p)}\ib  |\urm|^{\frac{p}{p-1}} \uhq \rho^{\frac{2-p}{p-1}}
 \Big)^{(p-1)/p} \Big( \epj\ib \frac{\uhq}{\rho^{2}} \Big)^{1/p} \overset{\footnotemark}{\leq}
 \frac{p}{p-1} \epj^{1/(1-p)}\ib  |\urm|^{\frac{p}{p-1}} \uhq \rho^{\frac{2-p}{p-1}}
+ \frac{\epj}{p}\ib \frac{\uhq}{\rho^{2}} .
\]

\footnotetext{Young: $ab\leq \frac{p-1}{p}a^{p/(p-1)}+
\frac{1}{p}b^{p}$.}

\no Applying H\"{o}lder inequality we get

\[
\ib  |\urm|^{\frac{p}{p-1}} \uhq \rho^{\frac{2-p}{p-1}} \overset{\footnotemark}{\leq} \Big( \ib |\urm|^{\frac{sp}{2(p-1)} } \rho^{\frac{(2-p)s}{2(p-1)}} \Big)^{2/s} \Big( \ib \uh^{\frac{qs}{s-2}} \Big)^{(s-2)/s}
\]
\footnotetext{H\"{o}lder: $(\frac{s}{2}, \frac{s}{s-2})$.}

\[
=
\Big( \ib |\urm|^{\frac{sp}{2(p-1)} } \rho^{\frac{(2-p)s}{2(p-1)}} \Big)^{2/s} \Big( \ib \uh^{q\frac{s-3}{s-2}}  \cdot \uh^{\frac{3q}{s-2}} \Big)^{(s-2)/s}
\]

\[
\overset{\footnotemark}{\leq} \Big( \ib |\urm|^{\frac{sp}{2(p-1)} } \rho^{\frac{(2-p)s}{2(p-1)}} \Big)^{2/s} \lqq{\uh}^{q\frac{s-3}{s}} \| \uh \|_{3q}^{q\frac{3}{s}}
\]

\footnotetext{H\"{o}lder: $(\frac{s-2}{s-3}, s-2)$.}
\[
\overset{\footnotemark}{\leq} \frac{3}{s} \epd \| \uh \|^{q}_{3q} + \frac{s-3}{s} \epd^{\frac{3}{3-s}} \Big( \ib |\urm|^{\frac{sp}{2(p-1)} } \rho^{\frac{(2-p)s}{2(p-1)}} \Big)^{2/(s-3)} \cdot \lqqq{\uh}.
\]

\footnotetext{Young: $(\frac{s}{3}, \frac{s}{s-3})$}

\no Thus we get

\[
\ib\frac{1}{\rho}\urm\uhq \leq  \frac{3(p-1)}{sp} \epj^{1/(1-p)}
 \epd \| \uh \|^{q}_{3q}
\]
\eqq{ +\frac{(p-1)(s-3)}{sp} \epj^{1/(1-p)}  \epd^{\frac{3}{3-s}}
\Big( \ib |\urm|^{\frac{sp}{2(p-1)} } \rho^{\frac{(2-p)s}{2(p-1)}}
\Big)^{2/(s-3)} \cdot \lqqq{\uh}+ \frac{\epj}{p}\ib
\frac{\uhq}{\rho^{2}} .}{u}

\no From Sobolev embedding theorem we have

\[
\| \uh \|^{q}_{3q} \leq c(q) \ib \pqp ,
\]

\no thus using (\ref{r}) and (\ref{u}) we get

\[
\pt  \lqqq{\uh}
 + \nu \frac{4(q-1)}{q}  \ib \pqp +\nu q \ib \frac{\uhq}{\rho^{2}}
\]
\[
 \leq  \frac{3(p-1)qc(q)}{sp} \epj^{1/(1-p)}
 \epd \ib \pqp + q\lqqq{\uh} + \lqqq{\hh}
\]
\[
+\frac{q(p-1)(s-3)}{sp} \epj^{1/(1-p)}  \epd^{\frac{3}{3-s}} \Big(
\ib |\urm|^{\frac{sp}{2(p-1)} } \rho^{\frac{(2-p)s}{2(p-1)}}
\Big)^{2/(s-3)} \cdot \lqqq{\uh}+ \frac{q\epj}{p}\ib
\frac{\uhq}{\rho^{2}} .
\]

\no Now we choose $\epj$ and $\epd$ small enough and we obtain

\[
\pt  \lqqq{\uh}
 + \nu \frac{2(q-1)}{q}  \ib \pqp +\nu q/2 \ib \frac{\uhq}{\rho^{2}}
\]
\[
 \leq \lqqq{\hh} + \Big[q +c \Big(
\ib |\urm|^{\frac{sp}{2(p-1)} } \rho^{\frac{(2-p)s}{2(p-1)}}
\Big)^{2/(s-3)}\Big] \lqqq{\uh},
 \]

\no where $c$ is a constant, which depends only on $q,p$ and $s$. If
we denote
\[
d(t):=\Big[q +c \Big( \ib |\urm|^{\frac{sp}{2(p-1)} }
\rho^{\frac{(2-p)s}{2(p-1)}} \Big)^{2/(s-3)}\Big],
\]

\no then we have \eqq{\pt  \lqqq{\uh}
 + \nu \frac{2(q-1)}{q}  \ib \pqp +\nu q/2 \ib \frac{\uhq}{\rho^{2}}
 \leq \lqqq{\hh} + d(t) \lqqq{\uh}.}{w}

\no  From (\ref{s}) and the assumption (\ref{m}) we know, that
function $d(t)$ in integrable on $ (t_{0}- \tau, t_{0})$. In
particular we have

\[
\pt  \lqqq{\uh}
 \leq \lqqq{\hh} + d(t) \lqqq{\uh}.
\]

\no If we multiply the sides by $\exp(-\ids)$ and integrate over
$(t_{0}-\tau, t)$, then we obtain the following estimate

\[
\| \uh(t) \|^{q}_{q} \leq e^{\ids} \| \uh (t_{0}- \tau)\|^{q}_{q} +
t \| h \|^{q}_{q, \infty}e^{\ids},
\]

\no i.e.

\eqq{\| \uh(t) \|_{q} \leq \m{const} \hd \m{ for }t\in (t_{0}- \tau,
t_{0}).}{x}

\no \textbf{Step 2.} We take $\ep \in (0,1)$ and we multiply the
sides of (\ref{j}) by $\frac{\oh}{\rde}$ and
integrate
\m{over $B_{2}$}

\[
\ib \oht \frac{\oh}{\rde} +\ib \ur \ohr \frac{\oh}{\rde}+ \ib \uz
\ohz\frac{\oh}{\rde} - \ib \frac{\ur}{\rte} \oh^{2} -2\ib
\frac{\uh}{\rho} \oor \frac{\oh}{\rde}
\]
\[
= \ib \gh \frac{\oh}{\rde}+ \nu  \ib \frac{1}{\rho} \poch{}{\rho}
(\rho \ohr) \frac{\oh}{\rde}+ \nu\ib \pochd{\oh}{z}\frac{\oh}{\rde}-
\nu\ib \frac{\oh^{2}}{\rce}.
\]

\no Now we can calculate

\[
\ib \oht \frac{\oh}{\rde}= \pohkrde,
\]
\[
\ib \ur \ohr \frac{\oh}{\rde}= \jd \ib \frac{\ur}{\rde} \ohkr= \jd
\ibt \frac{\ur}{\rje} \ohkr= - \jd \ibt (\frac{1}{ \rje}\urr+(\ep-1)
\frac{\ur}{\rde} ) \ohk
\]
\[
=- \jd \ib \frac{\ohk}{ \rde}\urr-(\frac{\ep}{2}-\jd) \ib
\frac{\ur}{\rho} \ohkrde ,
\]
\[
\ib \uz \ohz \frac{\oh}{\rde}=\jd \ib \frac{\uz}{\rde} \ohkz= -\jd
\ib \uzz \ohkrde ,
\]
\[
 \ib \frac{1}{\rho} \poch{}{\rho}
(\rho \ohr) \frac{\oh}{\rde}=  \ibt  \poch{}{\rho} (\rho \ohr)
\frac{\oh}{\rde}= - \ibt \frac{1}{\rje} \ohrk- \jd (\ep-2) \ibt
\ohkr \frac{1}{\rde}=
\]
\[
- \ibt \frac{1}{\rje} \ohrk +\jd (\ep-2)^{2} \ib \frac{\ohk}{\rce},
\]
\[
\ib \pochd{\oh}{z}\frac{\oh}{\rde}= - \ib \frac{1}{\rde}\ohzk.
\]

\no Then we get

\[
\pohkrde - \jd \ib \big[ \urr+ \frac{\ur}{\rho} + \uzz  \big]
\ohkrde - \epp \ib \frac{\ur}{\rho} \ohkrde - 2\ib \uhpr \oor
\frac{\oh}{\rde}
\]
\[
=\ib \gh \frac{\oh}{\rde} + \nu [\jd(2-\ep)^{2}-1] \ib
\frac{\ohk}{\rce} - \nu \ib \big[ \ohrk + \ohzk \big]
\frac{1}{\rde}.
\]

\no If we use (\ref{f}), then we have

\[
\pohkrde + \nu \ib \big[ \ohrk + \ohzk \big] \frac{1}{\rde} =2\ib
\uhpr \oor \frac{\oh}{\rde} +\epp \ib \frac{\ur}{\rho} \ohkrde
\]
\eqq{+\nu [\jd(2-\ep)^{2}-1] \ib \frac{\ohk}{\rce}+\ib \gh
\frac{\oh}{\rde}. }{z}

\no Now we shall estimate the right hand side. We recall that $\oor=
- \uhz$ and we get

\[
2\ib \uhpr \oor \frac{\oh}{\rde} = - 2\ib \uhpr \uhz
\frac{\oh}{\rde} = - \ib \uhkz \frac{\oh}{\rte}=  \ib
\frac{\uhk}{\rte} \ohz
\]
\[
\leq \frac{\nu}{2} \ib \frac{1}{\rde} \ohzk + \frac{1}{2\nu} \ib
\frac{\uh^{4}}{\rce}.
\]

\no We notice that

\[
\ib \ohrjerk \roe = \ib \big[ \frac{1}{\rje} \ohr +(\ep-1)
\frac{\oh}{\rde} \big]^{2} \roe
\]
\[
= \ib \frac{1}{\rde} \ohrk +
(\ep-1) \ibt \rdeo \ohkr + (\ep-1)^{2} \ib \frac{\ohk}{\rce}
\]
\[
= \ib \frac{1}{\rde} \ohrk - (\ep-1)(\ep-2) \ib  \frac{\ohk}{\rce} +
(\ep-1)^{2} \ib \frac{\ohk}{\rce}
\]
\[
= \ib \frac{1}{\rde} \ohrk + (\ep-1) \ib \frac{\ohk}{\rce}.
\]

\no If we notice that

\[
\ib \rdeo \ohzk = \ib \ohrjezk \roe,
\]
\no then  from (\ref{z}) we get

\[
\pohkrde + \nu \ib \Big\{  \ohrjerk + \ohrjezk \Big\} \roe
\]
\[
\leq \frac{\nu}{2} \ib \frac{1}{\rde} \ohzk + \frac{1}{2\nu} \ib
\frac{\uh^{4}}{\rce} +\epp \ib \frac{|\ur|}{\rho} \ohkrde
\]
\eqq{+\nu \frac{\ep}{2}(\ep-2) \ib \frac{\ohk}{\rce}+\ib |\gh|
\frac{|\oh|}{\rde}. }{aa}

\no Clearly we have

\[
\ib |\gh| \frac{|\oh|}{\rde} \leq \big\|
\frac{\gh}{\rho}\big\|_{6/5} \big\| \frac{\oh}{\rje}
\big\|_{6}\overset{\footnotemark}{\leq} c_{1} \big\| \nabla \big(
\frac{\oh}{\rje} \big) \big\|_{2} \leq c_{2}\Big(\ib \Big\{ \ohrjerk
+ \ohrjezk \Big\} \roe \Big)^{1/2}
\]
\[
\leq c_{3} + \frac{\nu}{4} \ib \Big\{  \ohrjerk + \ohrjezk \Big\}
\roe.
\]

\no If we use this estimate in (\ref{aa}), then we obtain

\[
\pohkrde + \frac{\nu}{4} \ib \Big\{  \ohrjerk + \ohrjezk \Big\} \roe
\]
\eqq{\leq  \frac{1}{2\nu} \ib \frac{\uh^{4}}{\rce} +\epp \ib
\frac{|\ur|}{\rho} \ohkrde +\nu \frac{\ep}{2}(\ep-2) \ib
\frac{\ohk}{\rce}+c_{3}. }{ab}

\no Finally, if we take the limit $\ep \rightarrow 0^{+}$, then we
have

\eqq{\pohkrd + \frac{\nu}{4} \ib  \ohrjrk + \ohrjzk \leq
\frac{1}{2\nu} \ib \frac{\uh^{4}}{\rho^{4}} +c_{3}.
 }{ac}

\no \textbf{Step 3. } We multiply (\ref{d}) by $\uhtrd$ and
integrate over $B_{2}$

\[
\ib \poch{\uh}{t}\uhtrd +\ib \ur \uhr \uhtrd+ \ib \uz \uhz \uhtrd+
\ib \frac{\uh^{4}}{\rho^{3}} \ur
\]
\[
 = \ib \hh \uhtrd + \nu \ib \frac{\uh^{3}}{\rho^{3}} \poch{}{\rho} (\rho \uhr)
 +\nu \ib
\pochd{\uh}{z}\uhtrd- \nu \ib \frac{\uh^{4}}{\rho^{4}}.
\]

\no Now we calculate

\[
\ib \poch{\uh}{t}\uhtrd= \jc \puhcrd,
\]
\[
\ib \ur \uhr \uhtrd = \jc \ibt \frac{\ur}{\rho} \poch{\uhc}{\rho} =
- \jc \ib \urr \uhcrd + \jc \ib \frac{\ur}{\rho} \uhcrd,
\]
\[
\ib \uz \uhz \uhtrd = \jc \ib \frac{\uz}{\rho^{2}} \poch{\uhc}{z} =
- \jc \ib \uzz \uhcrd,
\]
\[
\ib \frac{\uh^{3}}{\rho^{3}} \poch{}{\rho} (\rho \uhr) = \ibt
\frac{\uh^{3}}{\rho^{2}} \poch{}{\rho} (\rho \uhr) = - 3\ib
\frac{\uhk}{\rho^{2}} \uhrk + \ib \uhcrc,
\]
\[
\ib \pochd{\uh}{z}\uhtrd = - 3 \ib \uhdrd \uhzk.
\]

\no Thus we have

\[
\jc \puhcrd - \jc \ib \urr \uhcrd + \jc \ib \frac{\ur}{\rho} \uhcrd
- \jc \ib \uzz \uhcrd + \ib \frac{\ur \uhc }{\rho^{3}}
\]
\[
=\ib \hh \uhtrd - 3\nu \ib \frac{\uhk}{\rho^{2}} \uhrk + \nu \ib
\uhcrc - 3 \ib \uhdrd \uhzk,
\]

\no hence
\[
\jc \puhcrd - \jc \ib [\urr + \frac{\ur}{\rho} + \uzz ] \uhcrd +
\frac{3}{2} \ib \uhcurrt + 3 \nu \ib [\uhrk+ \uhzk] \uhdrd =\ib \hh
\uhtrd.
\]

\no If we use (\ref{f}) then we have

\eqq{\jc \puhcrd + \frac{3}{2} \ib \uhcurrt + 3 \nu \ib [\uhrk+
\uhzk] \uhdrd =\ib \hh \uhtrd .}{ad}

\no On the other hand we can write

\[
\ib \Big[\poch{}{z} \big( \frac{\uhk}{\rho} \big) \Big]^{2}= \ib
\rod \big(\uhkz \big)^{2}= 4 \ib \uhdrd \uhzk
\]
\[
\ib \uhkrrk = \ib \big[ 2 \frac{\uh}{\rho} \uhr- \uhdrd \big]^{2} =
4 \ib \uhdrd \uhrk - \ibt \frac{1}{\rho^{2}} \poch{\uhc}{\rho} + \ib
\uhcrc =4 \ib \uhdrd \uhrk - \ib \uhcrc.
\]

\no Thus using these equalities in (\ref{ad}) we get

\[
\jc \puhcrd + \frac{3}{4}\nu \ib \Big(\uhkrrk + \uhkrzk \Big)+
\frac{3}{4} \nu \ib \uhcrc = - \frac{3}{2} \ib \uhcurrt +\ib \hh
\uhtrd.
\]

\no From Young inequality we have

\[
\ib \hh \uhtrd = \ib \frac{\uh^{3}}{\rho^{3}} \cdot \rho h \leq
\frac{\nu}{4} \ib \uhcrc + c \ib \rho^{4} \hh^{4},
 \]

\no hence

\eqq{\jc \puhcrd + \frac{3}{4}\nu \ib \Big(\uhkrrk + \uhkrzk \Big)+
\frac{1}{2} \nu \ib \uhcrc \leq  \frac{3}{2} \ib \uhcurbrt +c. }{ae}

\begin{remark}
 In steps 2 and 3 we do not use the assumption on higher regularity
 of $\ur$.
\end{remark}

\no \textbf{Step 4. } We multiply (\ref{ae}) by $\frac{2}{\nu^{2}}$

\[
\frac{1}{2\nu^{2}} \puhcrd + \frac{3}{2\nu} \ib \Big(\uhkrrk +
\uhkrzk \Big)+ \frac{1}{\nu} \ib \uhcrc \leq  \frac{3}{\nu^{2}} \ib
\uhcurbrt +c
\]

\no We add this inequality to (\ref{ac}) and we obtain

\[
\frac{1}{2\nu^{2}} \puhcrd +\pohkrd + \frac{3}{2\nu} \ib
\Big(\uhkrrk + \uhkrzk \Big)+ \frac{\nu}{4} \ib  \ohrjrk + \ohrjzk
\]
\[
+\frac{1}{2\nu} \ib \uhcrc\leq \frac{3}{\nu^{2}} \ib \uhcurbrt +c
\]

\no Proceeding similarly as in \cite{pokornynestupa} we deduce that
$\| \bl{\omega} \|_{L^{2}}$ is integrable on $(t_{0}- \tau, t_{0})$,
thus it implies the boundedness of $\| D \bl{u}\|_{L^{2}}$ on
$(t_{0}- \tau, t_{0})$, therefore $(x_{0},t_{0})$ cannot be a
singular point for $\bl{u}$.

\begin{remark}
The Theorem can be proved also in the case $b=\infty$. Then we have
to assume that $\frac{3}{a}+\gamma<1$. In the proof we put
$p=\frac{2a}{2a-\delta a-3}$, $s=3+\delta a$, where $\delta$ is such
that $\frac{3}{a}+\gamma= 1-\delta$ and $\delta \in (0,
\frac{2a-3}{a})$.

\end{remark}

\end{document}